\documentclass[11pt]{amsart}
\usepackage{amsfonts,amssymb,amsthm,eucal}

\usepackage[margin=1.25in]{geometry}

\newtheorem{thm}{Theorem}
\newtheorem{cor}[thm]{Corollary}
\newtheorem{lemma}[thm]{Lemma}
\newtheorem{prop}[thm]{Proposition}

\newcommand{\R}{\mathbb{R}}

\newcommand{\N}{\mathbb{N}}

\newcommand{\C}{\mathbb{C}}

\newcommand{\inprod}[2]{\left\langle #1, #2 \right\rangle}

\newcommand{\abs}[1]{\left\vert #1 \right\vert}
\newcommand{\norm}[1]{\left\Vert #1 \right\Vert}

\newcommand{\dfn}[1]{{\bf #1}}
\newcommand{\eps}{\varepsilon}

\DeclareMathOperator{\vol}{vol}
\newcommand{\cln}{\colon}

\allowdisplaybreaks[3]
\numberwithin{thm}{section}
\numberwithin{equation}{section}

\author{Mark W.\ Meckes}

\thanks{Research supported in part by NSF grant DMS-0902203.}
\address{Case Western Reserve University, Cleveland, OH 44106, U.S.A.}
\email{mark.meckes@case.edu}


\title{Positive definite metric spaces}

\begin{document}

\begin{abstract}
  Magnitude is a numerical invariant of finite metric spaces, recently
  introduced by T.\ Leinster, which is analogous in precise senses to
  the cardinality of finite sets or the Euler characteristic of
  topological spaces. It has been extended to infinite metric spaces
  in several a priori distinct ways. This paper develops the theory of
  a class of metric spaces, positive definite metric spaces, for which
  magnitude is more tractable than in general. Positive definiteness
  is a generalization of the classical property of negative type for a
  metric space, which is known to hold for many interesting classes of
  spaces.  It is proved that all the proposed definitions of magnitude
  coincide for compact positive definite metric spaces and further
  results are proved about the behavior of magnitude as a function of
  such spaces.  Finally, some facts about the magnitude of compact
  subsets of $\ell_p^n$ for $p \le 2$ are proved, generalizing results
  of Leinster for $p=1,2$, using properties of these spaces which are
  somewhat stronger than positive definiteness.
\end{abstract}

\maketitle



\section{Introduction} \label{S:intro}

Magnitude is a canonical numerical invariant of finite metric spaces
recently introduced by Tom Leinster \cite{Leinster-cafe,Leinster},
motivated by considerations from category theory.  The same notion
appeared earlier, although it was not really developed, in connection
with measuring biodiversity \cite{SP}.  Magnitude is analogous in a
precise sense to the Euler characteristic of topological spaces or
partially ordered sets, and to the cardinality of finite sets, and it
may be interpreted as the \emph{effective number} of points of a
space. The definition of magnitude was extended to infinite metric
spaces in various ways in the papers
\cite{LW,Willerton-computer,Willerton-homogeneous,Leinster}. In this
setting, magnitude turns out to have close connections (some proved,
and some only conjectural at present) to classical invariants of
geometric measure theory and integral geometry, including Hausdorff
dimension and intrinsic volumes of convex bodies and Riemannian
manifolds. This paper is devoted to developing the theory of a
particular class of metric spaces, \emph{positive definite} metric
spaces, for which the theory of magnitude is more tractable than in
general. Examples of positive definite metric spaces include many
spaces of interest, including all subsets of $L_p$ when $1\le p \le 2$,
round spheres, and hyperbolic spaces.

Given a finite metric space $(A,d)$, its \dfn{similarity matrix} is
the matrix $\zeta_A \in \R^{A \times A}$ given by $\zeta_A(x,y) =
e^{-d(x,y)}$.  A \dfn{weighting} for $A$ is a vector $w \in \R^A$ such
that $\zeta_A w = \mathbf{1}$, the vector indexed by $A$ whose entries
are all $1$; i.e., $\sum_{y \in A} e^{-d(x,y)} w(y) = 1$ for every
$x\in A$.  If a weighting $w$ for $A$ exists, then the \dfn{magnitude}
of $A$ is defined to be $\abs{A} = \sum_{x\in A} w(x)$. (It is easy to
check that if multiple weightings for $A$ exist, they give the same
value for $\abs{A}$.)  The reader is referred to \cite{Leinster} for
the category-theoretic motivation of this definition, and to
\cite{Leinster, LW, Willerton-computer,Willerton-homogeneous} for
discussions of various intuitive interpretations of magnitude.

As a function of an arbitrary finite metric space, magnitude may
exhibit a number of pathological behaviors, the most obvious of which
is that it may be undefined. One simple condition that prevents this
unpleasant situation (as well as other pathologies; see \cite[Section
  2.4]{Leinster} for a number of relevant results) is if the
similarity matrix $\zeta_A$ is positive definite. In that case
$\zeta_A$ is in particular invertible, and $w = \zeta_A^{-1}
\mathbf{1}$ is a weighting for $A$. A finite metric space $A$ is
called \dfn{positive definite} (respectively, \dfn{positive
  semidefinite}) if $\zeta_A$ is a positive definite (positive
semidefinite) matrix.  Besides the fact that magnitude is always
defined, other nice properties of the class of positive definite
finite metric spaces include that magnitude is positive and monotone.
That is, if $A$ is a positive definite space and $\emptyset \neq B
\subseteq A$ (so that $B$ is positive definite as well), then $0 < 
\abs{B} \le \abs{A}$.

Three different approaches to extending the definition of magnitude to
infinite spaces were taken in
\cite{LW,Willerton-computer,Willerton-homogeneous,Leinster}. One of
the purposes of this paper is to show that these approaches are
essentially equivalent in the presence of an appropriate positive
definiteness assumption.  To that end, an arbitrary metric space $A$
is defined to be \dfn{positive definite} (respectively \dfn{positive
  semidefinite}) if each of its finite subsets is positive definite
(positive semidefinite) with their induced metrics. Other aims of this
paper are to investigate the regularity of magnitude as a function of
a positive definite metric space, and to clarify somewhat which metric
spaces are and are not positive definite.

As will be seen in Section \ref{S:examples}, a natural strengthening
of positive definiteness is equivalent to the classical property of
negative type for metric spaces. Although the terminology is more
recent, negative type was introduced and studied by Menger
\cite{Menger} and Schoenberg \cite{Schoenberg0,Schoenberg1}, and is
well-studied in the literature on metric embeddings; see
e.g.\ \cite{WW,DL}.  Thus the theory of magnitude naturally leads back
to this classical notion.

For clarity, a \dfn{metric space} $A = (A,d)$ here consists of a
nonempty set $A$ equipped with a \dfn{metric} $d \cln A\times A \to
[0,\infty)$ such that
\begin{itemize}
\item $d(x,y) = 0$ if and only if $x=y$,
\item $d(x,y) = d(y,x)$ for every $x,y \in A$, and
\item $d(x,y) \le d(x,z) + d(z,x)$ for every $x,y,z \in A$.
\end{itemize}
The category-theoretic motivation for the definition of magnitude is
based in part on the observation by Lawvere \cite{Lawvere} (which will
not be explained here) that a metric space is a particular instance of
an \emph{enriched category}. As pointed out to the author by
T.\ Leinster, of the properties of $d$ above only the triangle
inequality and the fact that $d(x,x) = 0$ for every $x$ are necessary
to Lawvere's observation (which moreover even allows infinite
distances); whereas some classical results used in this paper, for
example \cite[Theorem 1]{Schoenberg1}, require the symmetry property
of $d$ but not the triangle inequality. Attention will therefore be
restricted to the classical definition of a metric space as given
above.

\medskip

The rest of this paper is organized as follows.  The remainder of this
section establishes some additional notation and terminology. Section
\ref{S:magnitude} shows the equivalence, for compact positive definite
metric spaces, of several proposed definitions of magnitude, and
investigates continuity properties of magnitude as a function of the
metric space. Section \ref{S:examples} discusses sufficient conditions
for positive definiteness, in particular showing the connection with
negative type, and presents some counterexamples. Finally, Section
\ref{S:ell_p^n} generalizes some results of Leinster \cite{Leinster}
about the magnitude of subsets of Euclidean space $\ell_2^n$ and
taxicab space $\ell_1^n$ to $\ell_p^n$ spaces for $p < 2$, using
properties of those spaces which are stronger than positive
definiteness.


\section*{Notation, terminology, and conventions}

It will be useful to consider two general simple transformations of a
metric on a fixed set. If $A = (A, d)$ is a metric space, $t \in
(0,\infty)$, and \(\alpha \in (0,1]\), then $tA$ is shorthand for the
metric space $(A,td)$ and \(A^\alpha\) is shorthand for the metric
space $(A,d^\alpha)$.

For a metric space $A$, $M(A)$ denotes the space of finite signed
Borel measures on $A$.  Unless otherwise specified, a \dfn{measure}
will always refer to a finite signed Borel measure.  Denote further by
$M_+(A)$ the cone of positive measures on $A$, by $FM(A)$ the space of
finitely supported signed measures on $A$, and by $FM_+(A)$ the cone
of finitely supported positive measures on $A$.  The space $M(A)$ is
equipped with the norm $\norm{\mu} = \abs{\mu}(A)$, where $\abs{\mu}
\in M_+(A)$ is the total variation of $\mu$.

If $(X,d)$ is a metric space and $A,B\subseteq X$, the \dfn{Hausdorff
  distance} between $A$ and $B$ is
\[
d_H(A,B) = \max \Bigl\{ \sup_{a\in A}d(a,B),\ 
  \sup_{b\in B} d(b,A) \Bigr\};
\]
it is easy to verify that this defines a metric on the class of
compact subsets of $X$. If $A$ and $B$ are metric spaces then the
\dfn{Gromov--Hausdorff distance} between $A$ and $B$ is
\[
d_{GH}(A,B) = \inf d_H\bigl(\varphi(A),\psi(B)\bigr),
\]
where the infimum is over all metric spaces $X$ and isometric
embeddings $\varphi\cln A \hookrightarrow X$ and $\psi\cln
B\hookrightarrow X$. It is a nontrivial result that this defines a
metric on the family of isometry classes of compact metric spaces; see
\cite[Chapter 3]{Gromov}.

The precise normalizations used for Fourier transforms will not be
important here, but for concreteness, the Fourier transform of a
measure $\mu$ on $\R^n$ is defined as the function
\[
\widehat{\mu} (\omega) = \int_{\R^n} e^{-i2\pi\inprod{x}{\omega}} \ d\mu(x),
\]
where $\inprod{\cdot}{\cdot}$ is the standard inner product on $\R^n$,
and the Fourier transform of $f \in L_1(\R^n)$ is the function
\[
\widehat{f}(\omega) = \int_{\R^n} f(x) e^{-i2\pi\inprod{x}{\omega}} \ dx.
\]

For $0 < p < \infty$, $L_p$ will be used as shorthand for the vector
space $L_p[0,1]$ of equivalence classes (under almost everywhere
equality) of measurable functions $x:[0,1]\to \R$ such that
\[
\norm{x}_p = \left(\int_0^1 \abs{x(t)}^p \ dt\right)^{1/p} < \infty.
\]
As is well-known, $\norm{\cdot}_p$ defines a quasinorm on $L_p$ which
is only a norm when $p\ge 1$; it is less well-known that when $0 <
p < 1$, $d(x,y) = \norm{x-y}_p^p$ defines a metric on $L_p$.  Below,
$L_p$ will be equipped with the metric $d(x,y) =
\norm{x-y}_p^{\min\{1,p\}}$ unless otherwise specified.  An isometry
involving these spaces is understood as a \emph{metric}-preserving
function and not a \emph{quasinorm}-preserving function when $p<1$
(see the comments following Proposition
\ref{T:Levy-generalized}). Similarly, $\ell_p^n$ denotes $\R^n$
equipped with the metric $d(x,y) = \norm{x-y}_p^{\min\{1,p\}}$,
where \( \norm{x}_p = \bigl(\sum_{j=1}^n \abs{x_j}^p \bigr)^{1/p} \)
for $0<p<\infty$ and $\norm{x}_\infty = \max_{1\le j \le n}
\abs{x_j}$.


\section*{Acknowledgements}

I thank Tom Leinster for many interesting and enlightening
conversations about the subject of this paper and related ideas. I
also thank Carsten Sch\"utt for pointing out the results of \cite{KS},
and the web sites the $n$-Category Caf\'e\footnote{{\tt
    http://golem.ph.utexas.edu/category/}} and
MathOverflow\footnote{{\tt http://mathoverflow.net}}, which were
central in starting my involvement in this project.


\bigskip

\section{Magnitude of positive definite spaces} \label{S:magnitude}

Three different definitions for the magnitude of an infinite metric
space $A$ have been proposed in \cite{Leinster, LW,
  Willerton-computer, Willerton-homogeneous}. The first definition of
$\abs{A}$ is as the supremum of the magnitudes of finite subspaces of
$A$ (see \cite[Definition 3.1.1]{Leinster}). This is unsatisfactory in
general since magnitude is not monotone with respect to inclusion
among arbitrary finite metric spaces (see \cite[Example
2.1.7]{Leinster}), so that this definition is not consistent with the
original one when restricted to finite spaces.  However,
\cite[Corollary 2.4.4]{Leinster} shows that this is not the case among
finite positive definite metric spaces, making this definition
reasonable for compact positive definite metric spaces (the scope assumed
in \cite[Definition 3.1.1]{Leinster}).

The second approach, taken in \cite{LW} and \cite{Willerton-computer},
is to consider a sequence of finite subspaces $\{A_k\}$ such that
$\lim_{k\to \infty}A_k = A$ in the Hausdorff distance, and then define
$\abs{A} = \lim_{k\to \infty}\abs{A_k}$.  This is unsatisfactory since
it is not clear a priori whether this limit is independent of the
approximating subspaces $\{A_k\}$.

The third approach, taken in \cite{Willerton-homogeneous}, is to
generalize the original definition of magnitude using measures for
weightings. A \dfn{weight measure} for $(A,d)$ is a finite signed
measure $\mu \in M(A)$ such that $\int_A e^{-d(x,y)} \ d\mu(y) = 1$
for every $x \in A$. If $A$ possesses a weight measure $\mu$, then the
magnitude of $A$ may be defined as $\mu(A)$. If $A$ possesses multiple
weight measures, it is easy to check that they give the same value for
the magnitude; however it is not clear how generally weight measures
exist. (If $A$ is a compact homogeneous space then $A$ has a weight
measure; see \cite[Theorem 1]{Willerton-homogeneous}. Other sufficient
conditions follow from Lemma \ref{T:pos-weight} and Corollary
\ref{T:pos-weight-measure} below.)  In \cite{Willerton-homogeneous},
Willerton showed that the magnitudes of intervals, circles, and Cantor
sets, as computed via weight measures, coincide with their
magnitudes as computed in \cite{LW} using the second approach.

The results of this section show that all these approaches to defining
magnitude yield the same value of magnitude for compact positive
definite metric spaces, and also develop some continuity properties of
magnitude on such spaces. It will be convenient to take yet a fourth
approach to the definition of magnitude, in terms of a Rayleigh-like
quotient expression which already appears, in the finite case, in
\cite{Leinster}, and develop its relationships to the three approaches
described above.

\medskip

Given a compact metric space $(A,d)$, define a bilinear form $Z_A$ on $M(A)$
by 
\[
Z_A(\mu,\nu) = \int_A \int_A e^{-d(x,y)} \ d\mu(x) \ d\nu(y).
\]
By Fubini's theorem, $Z_A$ is symmetric. Observe that if $\mu$ is a
weight measure for $A$ and $\nu \in M(A)$, then $Z_A(\mu,\nu) =
\nu(A)$.

If $A$ is a compact positive definite metric space, the
\dfn{magnitude} $\abs{A}$ of $A$ is defined to be
\begin{equation}\label{E:pd-mag}
\abs{A} = \sup \left\{ \frac{\mu(A)^2}{Z_A(\mu,\mu)} \middle\vert
 \mu \in M(A),\ Z_A(\mu,\mu) \neq 0 \right\}.
\end{equation}
In \cite[Proposition 2.4.3]{Leinster} it is proved that
when $A$ is finite and positive definite, this coincides with the
earlier definition (cf.\ Theorem \ref{T:m-mag} below).  Observe that
an immediate consequence of this definition is that if $B\subseteq A$
then $\abs{B} \le \abs{A}$.

It will also be useful to consider the quantity
\begin{equation}\label{E:pd-div}
\abs{A}_+ = \sup \left\{\frac{\mu(A)^2}{Z_A(\mu,\mu)} \middle\vert
 \mu \in M(A)_+,\ \mu \neq 0 \right\}.
\end{equation}
Note that if $A$ is positive definite, then $Z_A(\mu,\mu) > 0$
whenever $\mu$ is a nonzero positive measure.  The quantity
$\abs{A}_+$ is called the \dfn{maximum diversity} of $A$ because of an
interpretation related to theoretical ecology (see
\cite{Leinster-entropy} and the discussion at the end of \cite[Section
  2.4]{Leinster}).  For any compact positive definite metric space
$A$, it is easy to check that $\abs{A}_+ \le
\exp(\operatorname{diam}(A))$.

Compactness is a useful and natural-seeming condition to assume in
this context.  However, it is not clear that it is necessarily the
most natural condition to use.  If $A$ is an infinite set in which
each distinct pair of points is separated by a distance $r>0$, then
$A$ is a \emph{noncompact} positive definite metric space, which can
nevertheless sensibly be assigned a finite magnitude $e^r$ using the
first definition proposed above (see \cite[Section 3.1]{Leinster}).
On the other hand, it is unknown at present whether the magnitude of a
compact positive definite metric space can be infinite.  As in
\cite[Section 3]{Leinster}, attention will nevertheless be restricted
here to compact spaces.

\medskip

The following lemma is central to the results of this section. Recall
that if $(A,d)$ is a metric space and $f\cln A \to \R$ is uniformly
continuous, the \dfn{modulus of continuity} of $f$ is the function
$\omega_f\cln (0,\infty) \to [0, \infty)$ defined by
\[
\omega_f(\eps) = \sup\bigl\{ \abs{f(x) - f(y)} \big\vert x,y \in
A,\ d(x,y) < \eps \bigr\}.
\]

\begin{lemma} \label{T:fin-approx-meas} Let $A$ and $B$ be compact
  subspaces of a metric space $X$ and let $\mu \in M(A)$. For any
  $\eps > d_H(A,B)$ there exists a $\nu \in M(B)$ such that \( \nu(B)
  = \mu(A) \), \( \norm{\nu} \le \norm{\mu}, \) and for any uniformly
  continuous $f\cln X \to \R$,
  \begin{equation*}
  \abs{\int_A f \ d\mu - \int_{B} f \ d\nu} \le \norm{\mu}
  \omega_f(\eps).
  \end{equation*}
  Moreover, if $\mu$ is positive then $\nu$ can be taken to be
  positive.
\end{lemma}

\begin{proof}
  Since $d_H(A,B) < \eps$, each point of $A$ is within distance $\eps$
  from some point of $B$.  Let $x_1,\dotsc,x_N \in B$ be the centers
  of open $\eps$-balls which cover $A$. Then the disjoint Borel sets
  $U_1 = B(x_1,\eps)$ and
  \[
  U_j = B(x_j,\eps) \setminus \bigcup_{k=1}^{j-1} B(x_k,\eps)
  \quad \text{for }j=2,\dotsc,N
  \]
  also cover $A$. Let $\nu = \sum_{j=1}^N \mu(U_j \cap A) \delta_{x_j}$.
  Then
  \begin{equation*}\begin{split}
      \abs{\int_A f\ d\mu - \int_{B} f \ d\nu}
      & = \abs{ \sum_{j=1}^N \int_{U_j} \bigl(f(x)-f(y_j)) \ d\mu(x)}
      \le \sum_{j=1}^N \int_{U_j} \abs{f(x)-f(y_j)} \ d\abs{\mu}(x) \\
      & \le \sum_{j=1}^N \omega_f(\eps) \abs{\mu}(U_j \cap A) 
      = \omega_f(\eps) \abs{\mu}(A) = \norm{\mu} \omega_f(\eps).
    \end{split}\end{equation*}

  Furthermore, 
  \[
  \nu(B) = \sum_{j=1}^N \mu(U_j \cap A) = \mu(A)
  \]
  and
  \[
  \norm{\nu} = \abs{\nu}(B) = \sum_{j=1}^N \abs{\mu(U_j \cap A)}
  \le \sum_{j=1}^N \abs{\mu}(U_j \cap A) = \abs{\mu}(A) = \norm{\mu}.
  \qedhere
  \]
\end{proof}

\medskip

\begin{lemma} \label{T:psd}
A compact metric space $A$ is positive semidefinite if and only if
$Z_A$ is a positive semidefinite bilinear form on $M(A)$.  If $Z_A$ is
positive definite then $A$ is positive definite.
\end{lemma}

\begin{proof}
  Recall that by definition $(A,d)$ is positive (semi)definite if all
  of its finite subspaces are positive (semi)definite.  The ``if''
  parts follow by applying the positive (semi)definite bilinear form
  $Z_A$ to finitely supported signed measures.

  Now suppose that $A$ is positive semidefinite and let $\mu \in M(A)$
  and $\eps > 0$. Apply Lemma \ref{T:fin-approx-meas} with $B$ a
  finite $\eps$-net in $A$ to obtain $\nu \in M(B)$ with
  $\norm{\nu}\le \norm{\mu}$ such that
  \[
  \abs{\int_A e^{-d(x,y)} \ d\mu(y) - \int_A e^{-d(x,y)} \ d\nu(y)} 
  \le \eps \norm{\mu}
  \]
  for each $x\in A$. From this it follows that
  \begin{equation*}\begin{split}
      \abs{Z_A (\mu,\mu) - Z_A(\nu,\nu)}
      & \le \abs{\int_A \int_A e^{-d(x,y)} \ d\mu(y) \ d\mu(x)
        - \int_A \int_A e^{-d(x,y)} \ d\nu(y) \ d\mu(x)}\\
      & \quad + \abs{\int_A \int_A e^{-d(x,y)} \ d\mu(x) \ d\nu(y)
        - \int_A \int_A e^{-d(x,y)} \ d\nu(x) \ d\nu(y)} \\
      & \le \eps \norm{\mu}^2 + \eps \norm{\mu} \norm{\nu}
      \le 2 \eps \norm{\mu}^2.
   \end{split}\end{equation*}
 Since $B$ is a positive semidefinite finite metric space,
 $Z_A(\nu,\nu) \ge 0$, which implies $Z_A(\mu,\mu) \ge -2 \eps
 \norm{\mu}^2$. Since $\eps > 0$ was arbitrary, $Z_A(\mu,\mu)\ge 0$.
\end{proof}

\medskip

The next result, which generalizes \cite[Proposition 2.4.3]{Leinster},
shows the agreement of the present definition \eqref{E:pd-mag} with
the measure-theoretic definition of magnitude used in
\cite{Willerton-homogeneous}, whenever both definitions can be
applied.

\begin{thm} \label{T:m-mag} Suppose $A$ is a compact positive definite
  metric space.  The supremum in \eqref{E:pd-mag} is achieved for a
  measure $\mu$ if and only if $\mu$ is a nonzero scalar multiple of a
  weight measure for $A$.  If $\mu$ is a weight measure for $A$ then
  $\abs{A} = \mu(A)$.
\end{thm}

\begin{proof}
  Suppose first that $\mu$ is a weight measure for $A$.  (If $\mu$
  achieves the supremum in \eqref{E:pd-mag}, then so does any nonzero
  scalar multiple of $\mu$ by homogeneity.) By Lemma \ref{T:psd},
  $Z_A$ is a positive semidefinite bilinear form on $M(A)$, and
  therefore satisfies the Cauchy--Schwarz inequality.  Thus if $\nu \in
  M(A)$, then
  \[
  \nu(A) = Z_A(\mu,\nu) \le \sqrt{Z_A(\mu,\mu) Z_A(\nu, \nu)} =
  \sqrt{\mu(A) Z_A(\nu, \nu)},
  \]
  with equality if $\nu = \mu$, and so
  \[
  \mu(A) = \frac{\mu(A)^2}{Z_A(\mu, \mu)} 
  = \sup \left\{ \frac{\nu(A)^2}{Z_A(\nu,\nu)} \middle\vert
    \nu \in M(A),\ Z_A(\nu,\nu) \neq 0 \right\}.
  \]

  Now suppose that $\mu$ achieves the supremum in \eqref{E:pd-mag} and
  let $\nu \in M(A)$ satisfy $\nu(A) = 0$. Then for any $t \in \R$,
  \[
  Z_A(\mu,\mu) \le  Z_A(\mu + t \nu,\mu + t \nu) 
  = Z_A(\mu, \mu) + 2 Z_A(\mu,\nu) t + Z_A(\nu,\nu) t^2.
  \]
  Since $Z_A(\nu,\nu) \ge 0$ by Lemma \ref{T:psd}, this implies that
  $Z_A(\mu,\nu) = 0$. Applying this in the case that $\nu = \delta_x
  - \delta_y$ for arbitrary $x,y \in A$ yields 
  \[
  \int_A e^{-d(x,z)} \ d\mu(z) = \int_A e^{-d(y,z)} \ d\mu(z),
  \]
  and thus $\mu$ is a scalar multiple of a weight measure for $A$.  
\end{proof}

\medskip

Theorem \ref{T:fin-approx} below shows that the present definition of
magnitude \eqref{E:pd-mag} coincides, for positive definite spaces,
with the first proposed definition discussed above.

\begin{thm}\label{T:fin-approx}
For any positive definite compact metric space \(A\),
\begin{equation}\label{E:fin-approx}
\abs{A} = \sup \left\{ \frac{\mu(A)^2}{Z_A(\mu,\mu)} \middle\vert
 \mu \in FM(A),\ \mu \neq 0 \right\}
 = \sup \bigl\{ \abs{B} \big\vert B \subseteq A \text{ is finite}
 \bigr\}.
\end{equation}
and
\begin{equation}\label{E:fin-approx-+}
\abs{A}_+ = \sup \left\{ \frac{\mu(A)^2}{Z_A(\mu,\mu)} \middle\vert
 \mu \in FM_+(A),\ \mu \neq 0 \right\}
 = \sup \left\{ \abs{B}_+ \middle\vert B \subseteq A \text{ is finite}
 \right\}.
\end{equation}
\end{thm}

\begin{proof}
  Observe first that when $(A,d)$ is positive definite and $\mu \in
  FM(A)$, it follows that $Z_A(\mu,\mu) = 0$ only for $\mu = 0$. The
  second equality in \eqref{E:fin-approx} follows from
  \cite[Proposition 2.4.3]{Leinster} (or Theorem \ref{T:m-mag} above),
  which shows that for finite positive definite spaces, the present
  definition \eqref{E:pd-mag} of magnitude agrees with the original
  definition. The second equality in \eqref{E:fin-approx-+} is
  immediate from \eqref{E:pd-div}.

In both \eqref{E:fin-approx} and \eqref{E:fin-approx-+} the first
quantity is by definition greater than or equal to the second quantity.
Let $\mu$ be a given measure on $A$ such that $Z_A(\mu,\mu) \neq 0$,
and let $\eps > 0$.  Apply Lemma \ref{T:fin-approx-meas} with $B$ an
$\eps$-net in $A$ to obtain $\nu \in FM(A)$, which is positive if
$\mu$ is positive, such that $\nu(A)=\mu(A)$, $\norm{\nu} \le
\norm{\mu}$ and
  \[
  \abs{\int_A f \ d\mu - \int_A f \ d\nu} \le \norm{\mu}
  \omega_f(\eps)
  \]
  for every continuous $f\cln A\to \R$. 

  Define $f_\mu\cln A \to \R$ by $f_\mu(x) = \int_A e^{-d(x,y)}
  \ d\mu(y)$, and define $f_\nu \cln A \to \R$ analogously. Then
  \[
  \abs{f_\nu(x) - f_\nu(y)} \le \norm{\nu} d(x,y),
  \]
  and
  \[
    \abs{f_\mu(x) - f_\nu(x)} =
    \abs{\int_A e^{-d(x,y)} \ d\mu(y) - \int_B e^{-d(x,y)} \ d\nu(y)}
    \le \norm{\mu} \eps.
  \]
  Consequently,
  \begin{align*}
    \abs{Z_A(\mu,\mu) - Z_B(\nu,\nu)}
    &= \abs{\int_A f_\mu \ d\mu - \int_A f_\nu \ d\nu} \\
    &\le \abs{\int_A f_\mu \ d\mu - \int_A f_\nu \ d\mu}
    + \abs{\int_A f_\nu \ d\mu - \int_A f_\nu \ d\nu} \\
    &\le \norm{\mu}^2 \eps + \norm{\mu}\norm{\nu} \eps 
    \le 2 \norm{\mu}^2 \eps. 
  \end{align*}
  Therefore
  \[
  \frac{\nu(A)^2}{Z_A(\nu,\nu)}
  = \frac{\mu(A)^2}{Z_A(\nu,\nu)}
  \ge \frac{\mu(A)^2}{Z_A(\mu,\mu) 
    + 2 \norm{\mu}^2 \eps}.
  \]
  Since $\eps > 0$ was arbitrary,
  \[
  \sup \left\{ \frac{\nu(A)^2}{Z_A(\nu,\nu)} \middle\vert
  \nu \in FM(A),\ \nu \neq 0\right\}
  \ge  \frac{\mu(A)^2}{Z_A(\mu,\mu)},
  \]
  and if $\mu$ is positive the same holds for the supremum over $FM_+(A)$.
\end{proof}

\medskip

In some circumstances magnitude can be expressed in terms of functions
instead of measures. A positive measure $\rho$ on a metric space $A$
is called a \dfn{good reference measure} if $\rho(U) > 0$ for every
nonempty open $U \subseteq A$.  For example, if $A\subseteq \R^n$ is
the closure of its interior as a subset of $\R^n$, then Lebesgue
measure restricted to $A$ is a good reference measure.  This may be
generalized naturally in at least two ways.  A metric space $A$ is
called \dfn{homogeneous} if its isometry group acts transitively on
the points of $A$.  It is well known that a locally compact
homogeneous metric space possesses an isometry-invariant (Haar)
measure, which is a good reference measure when restricted to any
subset of $A$ which is the closure of its interior.  On the other
hand, if $A$ is a metric space whose every nonempty open subset has
the same Hausdorff dimension $\delta$, then $\delta$-dimensional
Hausdorff measure is a good reference measure on $A$.

Given $h\in L_1(A,\rho)$, $h\rho$ denotes the signed measure on $A$
defined by $(h\rho)(S) = \int_S h \ d\rho$. The proof of the following
result is analogous to the proof of Theorem \ref{T:fin-approx}.

\begin{prop}\label{T:cts-mag}
For any positive definite compact metric space \(A\) with a good
reference measure $\rho$,
\[
\abs{A} = \sup \left\{ \frac{(h\rho)(A)^2}{Z_A(h\rho,h\rho)} \middle\vert
h \in L_1(A,\rho),\ Z_A(h\rho,h\rho) \neq 0 \right\}
\]
and
\[
\abs{A}_+ = \sup \left\{ \frac{(h\rho)(A)^2}{Z_A(h\rho,h\rho)} \middle\vert
h \in L_1(A,\rho),\ h \ge 0,\ h \text{ is not $\rho$-a.e.\ } 0 \right\}
\]
\end{prop}

\medskip

The agreement of \eqref{E:pd-mag} with the definition of magnitude as
the limit of magnitudes of an approximating sequence of subspaces, as
in \cite{LW,Willerton-computer}, will follow from the next result,
which is of independent interest.

\begin{thm}\label{T:lsc}
  The function $A \mapsto \abs{A}$ (with values in $[1,\infty]$) is
  lower semicontinuous with respect to Gromov--Hausdorff distance on
  the class of compact positive definite metric spaces.
\end{thm}

\begin{proof}
  Let $(A,d)$ be a positive definite metric space with $\abs{A} <
  \infty$ and let $\eps > 0$ be given. (The case where $\abs{A} =
  \infty$ is handled similarly.) Pick a signed measure $\mu$ on $A$
  with $Z_A(\mu,\mu) \neq 0$ such that
  \[
  \abs{A} \le \frac{\mu(A)^2}{Z_A(\mu,\mu)}(1 + \eps).
  \]
  Now let $B$ be any other positive definite metric space with
  $d_{GH}(A,B) > 0$.  Without loss of generality one may assume that
  $A,B\subseteq X$ for some metric space $X$, and $0 < d_H(A,B) \le 2
  d_{GH}(A,B)$. Let $\nu \in M(B)$ be as guaranteed by Lemma
  \ref{T:fin-approx-meas} with $2 d_H(A,B)$ in place of Lemma
  \ref{T:fin-approx-meas}'s $\eps$.

  Define $f_\mu\cln X \to \R$ by $f_\mu(x) = \int_A e^{-d(x,y)}
  \ d\mu(y)$, and define $f_\nu$ analogously. Then
  \[
  \abs{f_\nu(x) - f_\nu(y)} \le \norm{\nu} d(x,y),
  \]
  and by Lemma \ref{T:fin-approx-meas},
  \[
    \abs{f_\mu(x) - f_\nu(x)} =
    \abs{\int_A e^{-d(x,y)} \ d\mu(y) - \int_B e^{-d(x,y)} \ d\nu(y)}
    \le 2 \norm{\mu} d_H(A,B).
  \]
  Consequently,
  \begin{align*}
    \abs{Z_A(\mu,\mu) - Z_B(\nu,\nu)}
    &= \abs{\int_A f_\mu \ d\mu - \int_B f_\nu \ d\nu} \\
    &\le \abs{\int_A f_\mu \ d\mu - \int_A f_\nu \ d\mu}
    + \abs{\int_A f_\nu \ d\mu - \int_B f_\nu \ d\nu} \\
    &\le 2 \norm{\mu}^2 d_H(A,B) + 2 \norm{\mu}\norm{\nu} d_H(A,B) \\
    &\le 8 \norm{\mu}^2 d_{GH}(A,B).
  \end{align*}
  Therefore
  \begin{align*}
  \abs{B} &\ge \frac{\nu(B)^2}{Z_B(\nu,\nu)}
  \ge \frac{\mu(A)^2}{Z_A(\mu,\mu) 
    + 8 \norm{\mu}^2 d_{GH}(A,B)} \\
  &\ge \left(1 + \frac{8\norm{\mu}^2}{Z_A(\mu,\mu)} d_{GH}(A,B)\right)^{-1}
  \cdot \frac{\abs{A}}{1+\eps}.
  \end{align*}
  So if $d_{GH}(A,B) \le \frac{Z_A(\mu,\mu)}{8\norm{\mu}^2}\eps$, then
  $\abs{B} \ge (1+\eps)^{-2} \abs{A}$.
\end{proof}

\medskip

In general, $\abs{A}$ is not a continuous function of $A$. Examples
2.2.8 and 2.4.9 of \cite{Leinster} give an example, due to
S.\ Willerton, of a metric space $A$ such that $tA$ is positive
definite for each $t>0$ and $\lim_{t\to 0^+} \abs{tA} = 6/5$, whereas
$\lim_{t\to 0^+} tA = \{*\}$, which has magnitude $1$. It is an open
question whether $A \mapsto \abs{A}$ is continuous when restricted to
compact subsets of a fixed positive definite space.  The
\emph{asymptotic conjectures} of \cite{LW} (see also \cite[Conjecture
  3.5.10]{Leinster}) would imply in particular that magnitude is
continuous when restricted to compact convex subsets of $\ell_2^n$.

Theorem \ref{T:lsc} and the monotonicity of magnitude for positive
definite spaces immediately imply the following result, which shows
that the present definition \eqref{E:pd-mag} agrees with the
definition of magnitude in terms of an approximating sequence of
subspaces, and in particular shows that the latter definition is
independent of the subspaces chosen.

\begin{cor}\label{T:Hausdorff-limit}
  If $A$ is a compact positive definite metric space and $\{A_k\}$ is
  a sequence of compact subspaces of $A$ such that $\lim_{k \to
    \infty} d_H(A_k,A) = 0$, then $\abs{A} = \lim_{k \to \infty}
  \abs{A_k}$.
\end{cor}

\medskip

Theorems \ref{T:m-mag} and \ref{T:fin-approx} and Corollary
\ref{T:Hausdorff-limit} completely explain the agreement of
calculations of magnitudes using different definitions in \cite{LW}
and \cite{Willerton-homogeneous}, since all of the spaces involved are
positive definite by \cite[Proposition 2.4.13]{Leinster} and Theorem
\ref{T:pd-examples} below.

\medskip

The remaining results of this section develop some additional
properties of maximum diversity, which yield information about
magnitude for a particular class of spaces.  A compact positive
definite metric space $A$ is called \dfn{positively weighted} if
$\abs{A} = \abs{A}_+$. A finite positive definite space is positively
weighted if and only if its weighting has only nonnegative
components; several results in \cite[Section 2.4]{Leinster} give
sufficient conditions for this property and properties of the
magnitude of such finite spaces.

The following lemma gives several sufficient conditions for a compact
positive definite metric space to be positively weighted. It will be
seen in Corollary \ref{T:pos-weight-measure} that the first of these
conditions is also necessary.

\begin{lemma}\label{T:pos-weight}
  Let $(A,d)$ be a compact positive definite metric space.  Under any
  of the following conditions, $A$ is positively weighted.
  \begin{enumerate}
  \item \label{I:pos-wt-meas} There is a positive weight measure on
    $A$.
  \item $A$ is a homogeneous metric space.
  \item \label{I:fin-pos-wt} Every finite subset of $A$ has a
    weighting with only nonnegative components.
  \item There is an isometric embedding of $A$ into $\R$, with the
    standard metric on $\R$.
  \item For every $x,y,z\in A$, $d(x,y) \le \max\{ d(x,z), d(z,y) \}$
    (i.e., $A$ is an ultrametric space).
  \end{enumerate}
\end{lemma}

\begin{proof}
  \begin{enumerate}
  \item This follows from Theorem \ref{T:m-mag}.
  \item A compact homogeneous space has a nonnegative weight measure
    (see \cite[Theorem 1]{Willerton-homogeneous} and the comments
    following it).
  \item This follows from Theorem \ref{T:fin-approx}.
  \item By \cite[Theorem 4]{LW} or \cite[Proposition
    2.4.13]{Leinster}, every finite subset of $\R$ has a nonnegative
    weighting.
  \item By \cite[Proposition 2.4.18]{Leinster} (originally proved in
    \cite{POP}), every finite ultrametric space has a nonnegative
    weighting. \qedhere
\end{enumerate}
\end{proof}

\medskip

It should be noted that every subset of $\R$ and every ultrametric
space is positive definite (see \cite[Proposition 2.4.13]{Leinster}
and \cite{VN}, respectively; also \cite[Proposition 2.4.18]{Leinster}
and Theorem \ref{T:pd-examples} below). There exist homogeneous spaces
which are not positive definite (see \cite[Example
2.1.7]{Leinster}). Example 2.4.16 in \cite{Leinster} shows that when
$n\ge 2$, not all compact subsets of $\ell_1^n$ are positively
weighted, and numerical calculations in \cite{Willerton-computer} show
that not all compact subsets of $\ell_2^n$ are positively weighted,
although $\ell_1^n$ and $\ell_2^n$ are positive definite (see
\cite[Theorems 2.4.14 and 2.5.3]{Leinster}; also Theorem
\ref{T:pd-examples} below).

\begin{prop}\label{T:div-sup}
  If $A$ is a compact positive definite metric space, then the
  supremum in the definition \eqref{E:pd-div} of $\abs{A}_+$ is
  achieved by some $\mu \in M_+(A)$.
\end{prop}

\begin{proof}
  Denote by $P(A) = \{ \mu \in M_+(A) \mid \mu(A) = 1\}$ the space of
  probability measures on $A$. A well-known consequence of the
  Banach--Alaoglu theorem is that when $A$ is compact, $P(A)$ is
  compact with respect to the weak-$*$ topology inherited from the
  action of $M(A)$ as the dual of the Banach space $(C(A),
  \norm{\cdot}_\infty)$.  This topology on $P(A)$ is metrized by the
  Wasserstein distance
  \[
  d_W(\mu,\nu) = \sup \left\{ \int_A f \ d\mu - \int_A f \ d\nu
    \middle\vert f \cln A \to \R \text{ is $1$-Lipschitz} \right\}
  \]
 (see e.g.\ \cite[Corollary 6.13]{Villani}). By homogeneity,
  \begin{equation}\label{E:inf-P(A)}
  \abs{A}_+ = \sup_{ \mu \in P(A)}\frac{1}{Z_A(\mu,\mu)}.
  \end{equation}

  In the notation of the proof of Theorem \ref{T:lsc}, for
  $\mu, \nu \in P(A)$,
  \[
    \abs{Z_A(\mu,\mu) - Z_A(\nu,\nu)} \le \abs{\int_A f_\mu \ d\mu -
      \int_A f_\mu \ d\nu} + \abs{\int_A f_\nu \ d\mu - \int_A f_\nu \
      d\nu}
  \]
  since $\int_A f_\mu \ d\nu = \int_A f_\nu \ d\mu$ by Fubini's
  theorem.  As seen in the proof of Theorem \ref{T:lsc}, $f_\mu$ and
  $f_\nu$ are both $1$-Lipschitz (because $\norm{\mu} = \norm{\nu} =
  1$), and therefore
  \[
  \abs{Z_A(\mu,\mu) - Z_A(\nu,\nu)} \le 2 d_W(\mu,\nu).
  \]
  By continuity and compactness, the supremum in \eqref{E:inf-P(A)}
  is achieved by some $\mu \in P(A)$, which then also achieves
  the supremum in \eqref{E:pd-div}.
\end{proof}

\medskip

\begin{cor}\label{T:pos-weight-measure}
  If $A$ is a positively weighted compact positive definite metric
  space, then there is a positive weight measure on $A$.
\end{cor}

\begin{proof}
  By Proposition \ref{T:div-sup}, there is a $\mu \in M_+(A)$ which
  achieves the supremum in \eqref{E:pd-mag}. By Theorem \ref{T:m-mag},
  $\mu$ is, up to a scalar multiple, a weight measure for $A$.
\end{proof}

\medskip

\begin{prop}\label{T:pos-mag-cont}
The function \( A \mapsto \abs{A}_+ \) is continuous with respect to
the Gromov--Hausdorff distance on the class of compact positive
definite metric spaces.
\end{prop}

\begin{proof}
Lower semicontinuity follows exactly as in the proof of Theorem
\ref{T:lsc}.  The argument for upper semicontinuity proceeds along
similar lines.

Given $A$ and $B$ with $d_{GH}(A,B) > 0$, we may assume that $A,B
\subseteq X$ for some metric space $X$ and $0 < d_H(A,B) \le 2
d_{GH}(A,B)$. Suppose that $\mu \in M_+(B)$ satisfies
\[
\abs{B}_+ = \frac{\mu(B)}{Z_B(\mu,\mu)}
\]
(as guaranteed by Proposition \ref{T:div-sup}) and construct $\nu \in
M_+(A)$ as in Lemma \ref{T:fin-approx-meas}.  Proceeding analogously
to the proof of Theorem \ref{T:lsc}, one obtains that
\[
\frac{\mu(B)^2}{Z_B(\mu,\mu)} \le
\frac{\nu(A)^2}{Z_A(\nu,\nu) - 8\norm{\mu}^2 d_{GH}(A,B)}.
\]
Since $\mu$ is positive, $\norm{\mu} = \mu(B) = \nu(A)$, and so
\begin{align*}
\abs{B}_+ &\le \frac{\nu(A)^2}{Z_A(\nu,\nu) - 8 \nu(A)^2 d_{GH}(A,B)}
  \le \left(1 - \frac{8 \nu(A)^2}{Z_A(\nu,\nu)} d_{GH}(A,B)\right)^{-1}
  \abs{A}_+ \\
& \le \left(1 - 8\abs{A}_+ d_{GH}(A,B)\right)^{-1}
  \abs{A}_+.
\end{align*}
So if \( d_{GH}(A,B) \le \frac{\eps}{8\abs{A}_+} \) for $0 < \eps <
1$, then \( \abs{B}_+ \le (1-\eps)^{-1}\abs{A}_+ \).
\end{proof}

\medskip

\begin{cor}\label{T:pos-weights-cont}
  Magnitude is continuous with respect to the Gromov--Hausdorff
  distance on the class of positively weighted compact positive
  definite metric spaces.
\end{cor}

A particular case of Corollary \ref{T:pos-weights-cont} is that
magnitude is continuous with respect to Gromov--Hausdorff distance for
compact subsets of $\R$. In this setting, the slightly weaker result
of continuity with respect to Hausdorff distance also follows from an
exact integral formula for the magnitude of a compact subset of $\R$
given in \cite[Proposition 3.2.3]{Leinster}.


\bigskip

\section{Examples and counterexamples of positive definite spaces}
\label{S:examples}

This section is divided into two subsections. The first investigates
sufficient conditions for positive definiteness of a metric space, in
particular relating it to the classical property of negative type, and
gives a number of examples of positive definite metric spaces.  The
second subsection gives some examples of metric spaces which are not
positive definite, in particular demonstrating that some natural
operations on metric spaces do not preserve positive definiteness.


\bigskip

\subsection{Sufficient conditions for positive definiteness} ~

A function $f\cln E \to \C$ on a vector space $E$ is called
\dfn{positive definite} if, for every finite nonempty $A \subseteq E$,
the matrix $[f(x-y)]_{x,y\in A} \in \C^{A \times A}$ is positive
\emph{semidefinite}; $f$ is called \dfn{strictly positive definite} if
$[f(x-y)]_{x,y\in A}$ is positive definite. (This inconsistency in
terminology is unfortunately well-established.) Thus a
translation-invariant metric $d$ on $E$ is a positive definite
(respectively, positive semidefinite) metric if and only if $x \mapsto
e^{-d(x,0)}$ is a strictly positive definite (positive definite)
function.

The following classical result connects positive definiteness to
harmonic analysis (see e.g.\ \cite{Katznelson}).

\begin{prop}[Bochner's theorem]\label{T:Bochner}
A continuous function $f\cln \R^n \to \C$ is positive definite if and only
if $f=\widehat{\mu}$ for some positive measure $\mu$ on $\R^n$.
\end{prop}

Bochner's theorem does not consider \emph{strictly} positive definite
functions, and thus cannot directly identify positive definite metrics
on $\R^n$. A theory of strictly positive definite functions is
presented in \cite[Chapter 6]{Wendland}, which contains several
counterparts to Bochner's theorem, including the following.

\begin{prop}[{\cite[Theorem 6.11]{Wendland}}]\label{T:strict-Bochner}
  Suppose $f\in L_1(\R^n)$ is continuous. Then $f$ is strictly
  positive definite if and only if $f$ is bounded and $\widehat{f}$
  is nonnegative and not uniformly $0$.
\end{prop}

Rather than being applied directly here, Proposition
\ref{T:strict-Bochner} will be combined with classical results to
prove Theorem \ref{T:stability} below, which allows positive
definiteness of a metric space $A$ to be deduced from positive
\emph{semidefiniteness} of rescalings of $A$.

A metric space $A$ is \dfn{stably positive (semi)definite} if $tA$ is
positive (semi)definite for every $t > 0$. Since the definition of
magnitude implicitly involves an arbitrary choice of scale (as
discussed at the beginning of \cite[Section 2.2]{Leinster}), stable
positive definiteness is arguably a more natural condition on a metric
space than positive definiteness.  The space $A$ is of \dfn{negative
  type} if $A^{1/2}$ is isometric to a subset of a Hilbert space.
Spaces of negative type have been studied extensively in the theory of
embeddings of metric spaces; see e.g.\ \cite{WW,DL}, or \cite[Section
  2]{DW} for a concise recent survey. (The terminology \emph{negative
  type} stems from an alternative characterization of such spaces
which will not be needed here.)

The following result shows that the theory of magnitude leads
naturally to the classical notion of negative type, and that the
literature on negative type gives many examples of positive definite
metric spaces.  As mentioned above, it is also a useful tool for
upgrading positive semidefiniteness to positive definiteness, allowing
one to avoid using Proposition \ref{T:strict-Bochner} or related
results explicitly.

\begin{thm}\label{T:stability}
  The following are equivalent for a metric space $A$.
  \begin{enumerate}
  \item \label{I:stability-pos-def} $A$ is stably positive definite,
    and thus every compact subset of $A$ has a defined (possibly
    infinite) magnitude.
  \item \label{I:stability-semidef} $A$ is stably positive
    semidefinite.
  \item \label{I:stability-small-psd} There is a sequence $\{t_k > 0
    \mid k \in \N\}$ with $\lim_{k\to \infty} t_k = 0$ such that
    $t_kA$ is positive semidefinite for every $k$.
  \item \label{I:stability-neg-type} $A$ is of negative type.
  \end{enumerate}
\end{thm}

\begin{proof}
  The implications (\ref{I:stability-pos-def}) $\Rightarrow$
  (\ref{I:stability-semidef}) $\Rightarrow$
  (\ref{I:stability-small-psd}) are trivial.  The equivalence of
  (\ref{I:stability-semidef}) and (\ref{I:stability-neg-type}) was
  proved in \cite[Theorem 1]{Schoenberg1} (although that paper was
  written long before the terminology used here was introduced).  The
  proof will be completed by showing that
  (\ref{I:stability-small-psd}) $\Rightarrow$
  (\ref{I:stability-semidef}) and (\ref{I:stability-neg-type})
  $\Rightarrow$ (\ref{I:stability-pos-def}).

  \medskip

  Assume the condition in (\ref{I:stability-small-psd}). It suffices
  to assume that $A$ is finite. A finite metric space $A$ is positive
  semidefinite if and only if $\lambda_{\min}(\zeta_A) \ge 0$, where
  $\lambda_{\min}$ denotes the smallest eigenvalue of a symmetric
  matrix. Since the smallest eigenvalue of a symmetric matrix is a
  continuous function of the matrix entries (see
  e.g.\ \cite[p.\ 370]{HJ}), the set of $t$ such that $tA$ is positive
  semidefinite is a closed subset of $(0,\infty)$.

  If $A$ is a finite positive semidefinite metric space, then by the
  Schur product theorem \cite[Theorem 7.5.3]{HJ}, $\zeta_{nA}$ is a
  positive semidefinite matrix for every positive integer
  $n$. Therefore the condition in (\ref{I:stability-small-psd})
  implies that the set of $t$ such that $tA$ is positive semidefinite
  is a dense subset of $(0,\infty)$, and so $A$ is stably positive
  semidefinite.

  \medskip  

  Suppose finally that $(A,d)$ is of negative type. It suffices again
  to assume that $A$ is finite.  Then there is a function \( \varphi
  \cln A \to \ell_2^n \) for some $n$ such that $d(x,y) =
  \norm{\varphi(x)-\varphi(y)}_2^2$ for every $x,y \in A$, so in
  particular $\varphi$ is injective.  For $t>0$ define $f\cln \R^n \to
  \R$ by $f(x) = \exp\bigl(-t \norm{x}_2^2\bigr)$.  Proposition
  \ref{T:strict-Bochner} implies that $f$ is a strictly positive
  definite function, which means that
  \[
  \zeta_{tA} 
  = \bigl[\exp(-td(x,y))\bigr]_{x,y\in A}
  = \bigl[f \bigl(\varphi(x)-\varphi(y)\bigr)\bigr]_{x,y\in A}
  \]
  is a positive definite matrix. Thus $A$ is stably positive definite.
\end{proof}

\medskip

The equivalence of negative type with stable positive
\emph{semidefiniteness}, due to Schoenberg \cite{Schoenberg1}, is
well-known to experts on embeddings of metric spaces, and has been
generalized in various directions. The theory of magnitude, however,
requires positive \emph{definiteness}, and the equivalence of stable
positive definiteness with negative type appears to be new.  A further
equivalence, between negative type and Enflo's notion of generalized
roundness \cite{Enflo}, was proved in \cite{LTW}.

After Theorem \ref{T:stability} was first proved (with a slightly
weaker version of condition (\ref{I:stability-small-psd})),
T.\ Leinster found a direct proof of the implication
(\ref{I:stability-small-psd}) $\Rightarrow$
(\ref{I:stability-pos-def}).  This gives an alternative, more
elementary proof of the equivalence of conditions
(\ref{I:stability-pos-def})--(\ref{I:stability-small-psd}) which is
independent of the notion of negative type. Leinster's argument
furthermore obviates the need for Proposition \ref{T:strict-Bochner},
or indeed any mention of \emph{strictly} positive definite functions,
in the development of this theory. (However, the results of Section
\ref{S:ell_p^n} below require properties of certain functions which
are stronger than strict positive definiteness.) Leinster's proof is
included here with his permission.

\begin{proof}[Second proof of Theorem \ref{T:stability}, (\ref{I:stability-small-psd})
  $\Rightarrow$ (\ref{I:stability-pos-def})] Assume without loss of
  generality that $A$ is finite. By \cite[Proposition 2.2.6
    (i)]{Leinster}, $\zeta_{tA}$ is invertible for all but finitely
  many values of $t>0$, and thus there is an $\eps > 0$ such that
  $\zeta_{tA}$ is invertible for all $t \in (0,\eps)$. By the
  continuity property used above, either (a)
  $\lambda_{\min}(\zeta_{tA}) < 0$ for all $t \in (0,\eps)$ or (b)
  $\lambda_{\min} (\zeta_{tA}) > 0$ for all $t \in (0,\eps)$. By
  condition (\ref{I:stability-small-psd}), (a) is impossible, so (b)
  holds. Thus $tA$ is positive definite for all $t \in (0,\eps)$. By
  the positive \emph{definite} version of the Schur product theorem,
  $ntA$ is positive definite for every $t\in (0,\eps)$ and positive
  integer $n$, and thus $A$ is stably positive definite.
\end{proof}

\medskip

The remainder of this subsection is devoted to collecting examples of
metric spaces which are of negative type and hence positive definite.
The following result essentially goes back to L\'evy, generalized by
Bretagnolle, Dacunha-Castelle, and Krivine \cite{BDK}. For the precise
version stated here, see \cite[Theorem 8.9]{BL}.

\begin{prop}\label{T:Levy-generalized}
  Let $(E,\norm{\cdot}_E)$ be a real separable quasinormed space and
  let $0<p\le 2$.  The following are equivalent.
\begin{enumerate}
\item The function $f(x)=\exp(-\norm{x}_E^p)$ is positive definite on
  $E$.
\item \label{I:quasinorm-isometry} There is a linear map $T \cln  E \to
  L_p$ such that $\norm{x}_E = \norm{T(x)}_p$ for every $x \in E$.
\end{enumerate}
\end{prop}

\medskip 

Some remarks are in order at this point.  First, a map as in part
(\ref{I:quasinorm-isometry}) of Proposition \ref{T:Levy-generalized}
is usually called an isometry between quasinormed spaces, even when
the quasinorms are not norms. The slightly pedantic formulation of
Proposition \ref{T:Levy-generalized} above is to avoid ambiguous uses
of the word isometry.  Observe that if a quasinormed space
$(E,\norm{\cdot}_E)$ satisfies the conditions of Proposition
\ref{T:Levy-generalized} for some $p \in (0,1]$, then $d(x,y) =
\norm{x-y}_E^p$ defines a metric on $E$ and the map $T$ is an isometry
from $E$ equipped with this metric into $L_p$ with the metric $d(x,y)
= \norm{x-y}_p^p$.

Second, the restriction to separable spaces here is merely for
convenience of exposition (in order to avoid introducing nonseparable
$L_p$ spaces). The main interest here is whether compact subsets of
$E$ are positive definite, so this is no real restriction.

Finally, if a vector space $E$ has a metric $d$ which is homogeneous
of any degree (in the sense that for some $\beta > 0$, $d(tx,ty) =
t^\beta d(x,y)$ for every $x,y\in E$ and $t>0$), then $tE = (E,td)$ is
isometric to $E$ for every $t > 0$. It follows from Proposition
\ref{T:stability} that $(E,d)$ is stably positive definite if and only
if $(E,d)$ is positive semidefinite.

With these remarks in mind, the following is an immediate consequence
of Proposition \ref{T:Levy-generalized} and Theorem \ref{T:stability}.

\begin{cor}\label{T:pd-norm}
  The following are equivalent for a real separable normed space
  $(E,\norm{\cdot}_E)$.
  \begin{enumerate}
    \item $E$ is a positive semidefinite metric space.
    \item $E$ is a positive definite metric space.
    \item $E$ is a metric space of negative type.
    \item $E$ is isometric to a linear subspace of $L_1$.
  \end{enumerate}
\end{cor}

\medskip

The next result collects several large classes of metric spaces which
are known to have negative type, and whose compact subspaces therefore
all have well-defined magnitudes by Theorem \ref{T:stability}.  The
list is not intended to be exhaustive.

\begin{thm}\label{T:pd-examples}
  Any metric space from each of the following classes is of negative
  type, and hence positive definite.
  \begin{enumerate}
  \item \label{I:d^alpha} $A^\alpha$, where $A$ is a metric space
    of negative type and $0 < \alpha \le 1$.
  \item \label{I:Lp-pd} $L_p$ for $0 < p \le 2$ (with the metric
    $d(x,y) = \norm{x-y}_p^{\min\{1,p\}}$).
  \item \label{I:2d} Two-dimensional real normed spaces.
  \item \label{I:4-points} Metric spaces with at most four points.
  \item \label{I:ultrametric} Ultrametric spaces.
  \item \label{I:spheres} Round spheres (with the geodesic distance).
  \item \label{I:hyperbolic} Real or complex hyperbolic space.
  \item \label{I:trees} Weighted trees.
  \end{enumerate}
\end{thm}

\begin{proof}
\begin{enumerate}
\item This was proved by Schoenberg \cite[p.\ 527]{Schoenberg1}.

\item It was proved in \cite{BDK} that if $1\le p \le 2$, then $L_p$
  is isometric to a subspace of $L_1$.  Thus it suffices to prove the
  claim in the case $0 < p \le 1$.

  Suppose now that $0 < p \le 1$. The function $f(x) =
  \exp(-\norm{x}_p^p) = e^{-d(x,0)}$ is positive definite on $L_p$ by
  Proposition \ref{T:Levy-generalized}, which as observed above
  implies that $L_p$ is positive semidefinite and hence (by
  homogeneity and Theorem \ref{T:stability}) of negative type.

\item By \cite[Corollary 6.8]{Koldobsky}, every two-dimensional real
  normed space is isometric to a subspace of $L_1$.

\item By \cite[Theorem 1]{Wolfe}, every $N$-point metric space with $N
  \ge 4$ embeds isometrically into $\ell_\infty^{N-2}$.  In
  particular, every four-point space is isometric to a subset of
  $\ell_\infty^2$, which is two-dimensional so part (\ref{I:2d})
  applies, or more simply is isometric to $\ell_1^2$.

\item Every finite ultrametric space embeds isometrically in
  $\ell_2^n$ for some $n$; see e.g.\ \cite{Lemin}.
  
\item This follows from results in \cite{Kelly1,Kelly2}; see
  \cite[p.\ 263]{HLMT}.
 
\item This is proved in \cite[Corollaires 7.4 and 7.7]{FH}.

\item This is proved in \cite[Corollary 7.2]{HLMT}.
\qedhere
\end{enumerate}
\end{proof}

\medskip

The class of functions of metrics which preserve negative type, as
$t\mapsto t^\alpha$ does in part (\ref{I:d^alpha}) above, was
determined in \cite{Schoenberg2}.

Leinster proved directly that $\ell_1^n$ and $\ell_2^n$ are positive
definite for every $n$ in Propositions 2.4.14 and 2.5.3 of
\cite{Leinster}, respectively. Leinster's proof for $\ell_2^n$ is
based on the same idea behind the proof of Proposition
\ref{T:strict-Bochner} which partly underlies the (first) proof of
Theorem \ref{T:stability}.

The space $L_0$ of measurable functions $f \cln [0,1] \to \R$ is also
of negative type when equipped with an appropriate metric that
metrizes the topology of convergence in measure; see
\cite[p.\ 187]{BL}.

A direct proof that three-point metric spaces are positive definite is
given in \cite[Proposition 2.4.15]{Leinster}. No proof that every
four-point space has a defined magnitude is known which does not rely
on an embedding into a positive definite normed space.

Ultrametric spaces were directly proved to be positive definite in
\cite{VN}, see also \cite[Proposition 2.4.18]{Leinster}.

Some necessary conditions for manifolds to have negative type are also
known.  For example, a compact Riemannian manifold of dimension at
least two of negative type must be simply connected \cite[Theorem
  5.4]{HKM}, and a compact Riemannian symmetric space of negative type
must be a round sphere \cite[Corollary 2.6]{Kokkendorff}.  Since
compact symmetric spaces are homogeneous, their magnitude can
nevertheless be defined via weight measures as in
\cite{Willerton-homogeneous}.


\bigskip

\subsection{Counterexamples}
\label{S:nonexamples} ~

This subsection collects several examples which demonstrate the limits
of some of the results in the previous subsection, or show that some
appealing conjectures about positive definiteness are false.  Many of
these and related examples are known in the literature on metric
spaces of negative type.

The first example shows that the converse of Theorem
\ref{T:pd-examples} (\ref{I:d^alpha}) is false. Koldobsky
\cite{Koldobsky-1/2} constructed a normed space $E$ which embeds as a
quasinormed space (i.e.\ in the sense of Proposition
\ref{T:Levy-generalized} (\ref{I:quasinorm-isometry})) into $L_{1/2}$
but does not embed in $L_1$. By Proposition \ref{T:Levy-generalized}
and homogeneity, this means that the metric space $E^{1/2}$ is stably
positive definite, but $E$ is not positive semidefinite.

\medskip

The second example shows that the threshold $p=2$ in Theorem
\ref{T:pd-examples} (\ref{I:Lp-pd}) and the dimension two in Theorem
\ref{T:pd-examples} (\ref{I:2d}) are both optimal.  It was proved by
Dor \cite{Dor} that if $2 < p \le \infty$ then $\ell_p^3$ does not
embed isometrically in $L_1$, and is therefore not positive definite.
In particular, $L_p$ is not positive definite for any $p>2$.

\medskip

The third example shows that the cardinality four in Theorem
\ref{T:pd-examples} (\ref{I:4-points}) is optimal. As shown in
\cite[Example 2.2.7]{Leinster}, the vertices of the complete bipartite
graph $K_{3,2}$, with the shortest path metric, form a metric space
which is not positive definite if all edges have equal lengths $r <
\log \sqrt{2}$. 

\medskip

The fourth set of examples are compact Riemannian manifolds which are
not positive definite with the geodesic metric.  By \cite[Theorem
  5.4]{HKM}, if $M$ is any non-simply connected compact Riemannian
manifold $M$ of dimension at least two, then $M$ fails to have
negative type. Thus by Theorem \ref{T:stability} there is some $t>0$
such that $tM$ is not positive definite.

It is also possible to give an example which is both more elementary
and is topologically a sphere.  The idea is simply to construct a
surface $S$ which almost isometrically contains a copy of $K_{3,2}$
with short edge lengths.

Start with a $2$-sphere with radius smaller than $(\log 2)/\pi$ in
$\R^3$, and consider the following five points: two opposite poles,
and three equidistant points on the equator.  Drawing the lines of
longitude through the latter three points, we obtain a copy of
$K_{3,2}$ in the sphere with equal edge lengths $r < \log
\sqrt{2}$. Now put large bulges on the sphere in the three regions
delineated by the three lines of longitude.  These bulges may be made
large enough to make the geodesic distance on $S$ between the three
equatorial points arbitrarily close to the graph distance $2r$.

Since $K_{3,2}$ with edge lengths $r$ is not positive semidefinite,
and the smallest eigenvalue of a symmetric matrix is a continuous
function of the matrix entries, it follows that the bilinear form
$Z_S$ on $M(S)$ is not positive semidefinite.

\medskip

The class of positive definite metric spaces is closed under taking
$\ell_1$ products \cite[Lemma 2.4.2 (ii)]{Leinster}. The next set of
examples shows that $\ell_p$ products for any $p>1$ fail to preserve
positive definiteness, even in the more restricted context of positive
definite normed spaces (in which one usually speaks of $\ell_p$
\emph{sums}).  For $p>2$, this follows from the fact discussed above
that $\ell_p^3$, which is the $\ell_p$ product of three copies of the
positive definite space $\R$, is not positive definite. For $1<p\le
2$, \cite[Corollary 3.4]{KS} shows that the $\ell_p$ sum of $n$ copies
of $\ell_1^n$ does not embed isometrically in $L_1$ when $n$ is
sufficiently large (specifically, when $n > (5\sqrt{2})^{p/(p-1)}$,
although this bound is not sharp), and hence is not positive
semidefinite.

For a small concrete example, let $A = \{ 0, \pm e_1, \pm e_2\}
\subseteq \ell_1^2$, and consider $A \times A \subseteq \ell_1^2
\oplus_2 \ell_1^2$. Then a numerical calculation of the eigenvalues of
$\zeta_{t(A\times A)}$ for small $t$ shows that $A \times A$ is not
stably positive definite.

\medskip

The last examples show that the property of positive definiteness has
no simple relationship with $1$-Lipschitz maps. First, positive
definiteness is not preserved by $1$-Lipschitz maps. In fact, since
every separable Banach space is isometric to a quotient space of the
positive definite space $\ell_1$ (cf.\ \cite[p.\ 108]{LT}), it happens
quite generically that 1-Lipschitz maps fail to preserve positive
definiteness.  For a concrete low-dimensional example, one can define
a linear such map $\ell_1^4 \to \ell_\infty^3$ (recall that the former
space is positive definite and the latter is not) by mapping the
standard basis vectors of $\ell_1^4$ to the four vertices on one facet
of the $\ell_\infty^3$ unit ball.

On the other hand, if $B$ is a positive definite metric space and
there is a $1$-Lipschitz surjection $f \cln A \to B$, $A$ need not be
positive definite.  In the language of \cite[Definition
  2.2.4]{Leinster}, this means that positive definiteness is not
preserved by expansions.  One can even insist that $f$ be a bijection.
Let $E$ be a normed space which is not positive definite and let $A
\subseteq E$ be any finite subset which is not positive definite. Then
a generic linear functional $f \cln E \to \R$ is injective when
restricted to $A$, and $B = f(A) \subseteq \R$ is positive definite.


\bigskip

\section{Magnitude in $\ell_p^n$} \label{S:ell_p^n}

This section generalizes some results of Leinster \cite{Leinster}
about the magnitude of subsets of $\ell_1^n$ and $\ell_2^n$ to the
spaces $\bigl(\ell_p^n\bigr)^\alpha$ for $0 < p \le 2$ and $0 < \alpha
\le 1$. Leinster's proofs for Euclidean space are based on the exact
formula \eqref{E:Poisson-kernel} for the Fourier transform of the
function $x \mapsto e^{-\norm{x}_2}$ on $\R^n$ (for $\ell_1^n$ more
elementary tools suffice).  For $0 < p \le 2$, $0 < \alpha \le 1$, and
a positive integer $n$, define
\[
F_{p,\alpha}^n\cln \R^n \to \R, \qquad
F_{p,\alpha}^n(x) = \exp \bigl(-\norm{x}_p^{\alpha \min\{1,p\}} \bigr).
\]
The generalizations here require proving that the Fourier transforms
$\widehat{F_{p,\alpha}^n}$ share the properties of
$\widehat{F_{2,1}^n}$ which are essential in Leinster's proofs.

Parts (\ref{I:d^alpha}) and (\ref{I:Lp-pd}) of Theorem
\ref{T:pd-examples} imply that $\bigl(\ell_p^n\bigr)^\alpha$ is of
negative type, which by Theorem \ref{T:stability} is equivalent to the
statement that $F_{p,\alpha}^n$ is a strictly positive definite
function.  Lemma \ref{T:l_p-poly-bound} is a quantitative sharpening
of this fact. In probabilistic terms, it gives polynomial lower bounds
on the densities of a particular class of stable random vectors. As
such, it may already be known in the probability literature, although
we have been unable to find a statement. The bounds given by the proof
are nonoptimal, but are sufficient for the purposes of the present
paper.

\begin{lemma}\label{T:l_p-poly-bound}
Given $0 < p \le 2$, $0 < \alpha \le 1$, and a positive integer $n$ there
is a constant $c_{p,\alpha,n} > 0$ such that
\[
\widehat{F_{p,\alpha}^n}(\omega) \ge c_{p,\alpha,n} \bigl( 1 +
\norm{\omega}_2 \bigr)^{-(1+p)n}
\]
for every $\omega \in \R^n$.
\end{lemma}

\begin{proof}
Suppose first that $0 < p < 2$. Define $\gamma_p\cln \R \to \R$ by
$\gamma_p(t) = e^{-\abs{t}^p}$.  Then \cite[Lemma 2.27]{Koldobsky}
shows that $\widehat{\gamma_p} > 0$ everywhere, and asymptotic
expansions in Theorems 2.4.1, 2.4.2, and 2.4.3 in \cite{IL} imply that
\[
\lim_{\abs{\omega} \to \infty}
\widehat{\gamma_p}(\omega)\abs{\omega}^{1+p}
\]
exists and is finite.  It follows that there is a constant $c(p)>0$
such that 
\begin{equation}\label{E:gamma-bound}
\widehat{\gamma_p}(\omega) \ge c(p)
(1+\abs{\omega})^{-(1+p)}
\end{equation}
for every $\omega \in \R$.

A theorem of Bernstein \cite[Theorem
  XIII.4]{Feller} implies that for every $r \in (0,1]$ there is a
probability measure $\mu_r$ on $[0,\infty)$ such that
\begin{equation}\label{E:Bernstein}
e^{-t^r} = \int_0^\infty e^{-ts} \ d\mu_r(s)
\end{equation}
for every $t \ge 0$. In particular, for $r = \alpha \min \{1/p, 1\}$,
\[
F_{p,\alpha}^n (x) = \int_0^\infty e^{-s\norm{x}_p^p} \ d\mu_{r}(s)
= \int_0^\infty \left(\prod_{j=1}^n \gamma_p \bigl(s^{1/p} x_j \bigr) \right)
  \ d\mu_{r}(s).
\]

Now by Fubini's theorem, a linear change of variables, and
\eqref{E:gamma-bound},
\begin{align}
\widehat{F_{p,\alpha}^n}(\omega) &= \int_0^\infty
\left( \prod_{j=1}^n s^{-1/p} \widehat{\gamma_p}\bigl (s^{-1/p}\omega_j \bigr) \right)
\ d\mu_{r}(s) \label{E:Fpqd=} \\
&\ge c(p)^n \int_0^\infty
\left(s^{-1/p} \bigl(1+s^{-1/p}\norm{\omega}_2\bigr)^{-(1+p)}\right)^n
\ d\mu_{r}(s). \nonumber
\end{align}
If $0\le s \le 1$, then
\begin{equation*}
s^{-1/p} \bigl(1+s^{-1/p}\norm{\omega}_2\bigr)^{-(1+p)}
= s \bigl(s^{1/p} + \norm{\omega}_2\bigr)^{-(1+p)}
\ge s \bigl(1 + \norm{\omega}_2\bigr)^{-(1+p)}.
\end{equation*}
If $s>1$, then
\begin{equation*}
s^{-1/p} \bigl(1+s^{-1/p}\norm{\omega}_2\bigr)^{-(1+p)}
\ge s^{-1/p} \bigl(1+\norm{\omega}_2\bigr)^{-(1+p)}.
\end{equation*}
Thus
\[
\widehat{F_{p,q}^n(\omega)} \ge c(p)^n \left(\int_{[0,1]} s^n \ d\mu_{r}(s)
  + \int_{(1,\infty)} s^{-n/p} \ d\mu_{r}(s) \right)
\bigl(1+\norm{\omega}_2)^{-(1+p)d}.
\]

Now suppose $p=2$. When also $\alpha = 1$, there is the exact formula
\begin{equation}\label{E:Poisson-kernel}
\widehat{F_{2,1}^n}(\omega) = \frac{c_n}{\bigl(c_n' +
  \norm{\omega}_2^2\bigr)^{(n+1)/2}},
\end{equation}
where $c_n,c_n'>0$ are constants (which can be given explicitly)
depending only on $n$ (see \cite[Theorem I.1.4]{SW}).  By
\eqref{E:Bernstein}, Fubini's theorem, a linear change of variables,
and \eqref{E:Poisson-kernel},
\begin{equation}\label{E:F2=}
\widehat{F_{2,\alpha}^n}(\omega) = \int_0^\infty s^{-n} 
\frac{c_n}{\bigl(c_n' + s^{-2} \norm{\omega}_2^2\bigr)^{(n+1)/2}}
\ d\mu_\alpha(s).
\end{equation}
The proof is completed as before.
\end{proof}

\medskip

\begin{lemma}\label{T:monotonicity}
Given $0 < p \le 2$, $0 < \alpha \le 1$, and a positive integer $n$,
for each $\omega \in \R^n$, the function $t \mapsto
\widehat{F_{p,\alpha}^n}(t\omega)$ is decreasing for $t\ge 0$.
\end{lemma}

\begin{proof}
For $p=2$ the lemma follows from \eqref{E:F2=} in the proof of Lemma
\ref{T:l_p-poly-bound}.

For $p<2$, \cite[Theorem 2.5.3]{IL} implies that for any $\omega \in
\R$, $t \mapsto \widehat{\gamma_p}(t \omega)$ is decreasing on $[0,
  \infty)$.  The lemma in this case follows from this fact and the
  equality in \eqref{E:Fpqd=}.
\end{proof}

\medskip

The main results of this section, Theorems \ref{T:l_p-finite} and
\ref{T:growth-upper}, were proved by Leinster in the cases that
$\alpha = 1$ and $p=1,2$ (see Theorem 3.4.8, Proposition 3.5.3, and
Theorem 3.5.5 in \cite{Leinster}). The proofs below
generalize Leinster's proofs for $\ell_2^n$, using Lemmas
\ref{T:l_p-poly-bound} and \ref{T:monotonicity} in place of the exact
formula \eqref{E:Poisson-kernel}. The results of Section
\ref{S:magnitude} allow the exposition to be simplified somewhat by
working with measures instead of finite subsets.

\begin{thm}\label{T:l_p-finite}
Let $A$ be a compact subset of $\bigl(\ell_p^n\bigr)^\alpha$, where
$0<p\le 2$ and $0 < \alpha \le 1$. Then $\abs{A} < \infty$.
\end{thm}

\begin{proof}
  Let $\psi\cln \R^n \to \R$ be an even, compactly supported,
  $C^\infty$ function such that $\psi(x) = 1$ for all $x \in A-A = \{
  y - z \mid y,z\in A\}$. Then $\widehat{\psi}$ is a real-valued
  Schwartz function. By Lemma \ref{T:l_p-poly-bound} there is some
  constant $C(p,\alpha,\psi) > 0$ such that
\[
\widehat{\psi} \le C(p,\alpha,\psi) \widehat{F_{p,\alpha}^n}
\]
everywhere on $\R^n$.

Since $F_{p,\alpha}^n$ is positive definite, integrable, and
continuous at $0$, \cite[Corollary I.1.26]{SW} implies that
$\widehat{F_{p,\alpha}^n} \in L_1(\R^n)$ and thus (since
$\widehat{F_{p,\alpha}^n}$ is also even) that $F_{p,\alpha}^n$ is the
Fourier transform of $\widehat{F_{p,\alpha}^n}$. Now for any $\mu \in
M(A)$, by Fubini's theorem,
\begin{align*}
Z_A(\mu,\mu) & = \int_A \int_A F_{p,\alpha}^n (x-y) \ d\mu(x) \ d\mu(y) \\
&= \int_A \int_A \int_{\R^n} \widehat{F_{p,\alpha}^n}(\omega) 
  e^{-i2\pi\inprod{x-y}{\omega}} \ d\omega \ d\mu(x) \ d\mu(y) \\
&= \int_{\R^n} \abs{\widehat{\mu}(\omega)}^2
  \widehat{F_{p,\alpha}^n}(\omega) \ d\omega \\
& \ge \frac{1}{C(p,\alpha,\psi)} 
  \int_{\R^n} \abs{\widehat{\mu}(\omega)}^2
  \widehat{\psi}(\omega) \ d\omega \\
&= \frac{1}{C(p,\alpha,\psi)} \int_A \int_A \psi(x-y)
  \ d\mu(x) \ d\mu(y) = \frac{1}{C(p,\alpha,\psi)} \mu(A)^2.
\end{align*}
Thus by \eqref{E:pd-mag}, $\abs{A} \le C(p,\alpha,\psi)$.
\end{proof}

\medskip

It is natural, especially in light of Theorem \ref{T:pd-examples}, to
ask whether Theorem \ref{T:l_p-finite} applies to $(E)^\alpha$ for
arbitrary finite dimensional subspaces $E$ of $L_p$. To extend the
present proof to this setting would require generalizing a lower bound
as in Lemma \ref{T:l_p-poly-bound} to the densities of much more
general stable random vectors. Although such bounds have been the
subject of much study (see e.g.\ \cite{Watanabe}) it appears that
known results do not suffice for this purpose.

\medskip

Given a vector space $V$, we define
\[
t \cdot A = \{ ta \mid a \in A\}
\]
for each $t > 0$ and $A \subseteq V$. Observe that when $V$ is
equipped with a metric which is not homogeneous of degree $1$, the
abstract metric space $t A$ is not in general isometric to $t \cdot A
\subseteq V$.  (\emph{Homogeneous} is used here in the sense that for
some $\beta > 0$, $d(tx,0) = t^\beta d(x,0)$ for every $x\in \R^n$ and
$t>0$, and not in the sense of possessing a transitive isometry
group.)  In particular, if the metric on $V$ is homogeneous of degree
$\beta$, then $t A$ is isometric to $t^{1/\beta}\cdot A$.  This
includes in particular $V = \bigl(\ell_p^n\bigr)^\alpha$ in the entire
range $0<p\le 2$ and $0 < \alpha \le 1$, for which $\beta = \alpha
\min\{1,p\}$.

\begin{thm}\label{T:growth-upper}
  Let $A$ be a compact subset of $\bigl(\ell_p^n\bigr)^\alpha$, where
  $0<p\le 2$ and $0 < \alpha \le 1$. There is a constant $C > 0$ such
  that $\abs{tA} \le C t^{n/\beta}$ for all $t\ge 1$, where $\beta =
  \alpha \min\{1,p\}$.
\end{thm}

\begin{proof}
Let $\psi$ be as in the proof of Theorem \ref{T:l_p-finite}, and for
$t\ge 1$ define $\psi_t(x) = \psi(\frac{x}{t})$, so that $\psi_t$ is
an even, compactly supported, nonnegative $C^\infty$ function such
that $\psi_t(x) = 1$ for all $x \in \{ ty - tz \mid y,z\in A\}
$. Lemma \ref{T:monotonicity} and the proof of Theorem
\ref{T:l_p-finite} imply that
\[
\abs{t\cdot A}
\le \sup_{\omega \in \R^n} \frac{\widehat{\psi_t}(\omega)}
  {\widehat{F_{p,\alpha}^n}(\omega)}
= \sup_{\omega \in \R^n} \frac{t^n \widehat{\psi}(t\omega)}
  {\widehat{F_{p,\alpha}^n}(\omega)}
= t^n \sup_{\omega \in \R^n} \frac{\widehat{\psi}(\omega)}
  {\widehat{F_{p,\alpha}^n}(\omega/t)}
\le t^n \sup_{\omega \in \R^n} \frac{\widehat{\psi}(\omega)}
  {\widehat{F_{p,\alpha}^n}(\omega)}. 
\]
The theorem follows by replacing $t$ with $t^{1/\beta}$.
\end{proof}

\medskip

The last theorem of this section is a generalization of another result
of Leinster \cite[Theorem 3.5.6]{Leinster} which complements Theorem
\ref{T:growth-upper} for subsets with positive volume. Leinster states
the result for finite dimensional positive definite normed spaces, but
the proof, which will not be repeated here, generalizes immediately
from norms to translation invariant, homogeneous metrics on
$\R^n$.
\begin{thm}\label{T:growth-lower}
Let $d$ be a positive definite, translation invariant metric on $\R^n$
which is homogeneous of degree $\beta \in (0,1]$, and let $B = \{ x
\in \R^n \mid d(x,0) \le 1\}$. If $A \subseteq (\R^n, d)$ is compact,
then
\[
\abs{A} \ge \frac{\vol(A)}{\Gamma\bigl(\frac{n}{\beta} + 1\bigr)
  \vol(B)}.
\]
In particular,
\[
\abs{t A} = \abs{t^{1/\beta} \cdot A} \ge
\frac{\vol(A)}{\Gamma\bigl(\frac{n}{\beta} + 1\bigr) \vol(B)}
t^{n/\beta}
\]
for every $t>0$.
\end{thm}

\medskip

In the language of \cite[Definition 3.4.5]{Leinster}, Theorem
\ref{T:growth-upper} implies that any compact subset $A\subseteq
\bigl(\ell_p^n\bigr)^\alpha$ has magnitude dimension at most
$n/\beta$, and Theorem \ref{T:growth-lower} implies that if $A$ also
has positive volume, then the magnitude dimension of $A$ is precisely
$n/\beta$.\footnote{The published version of this paper used the
  notation $tA$ in the last two results to mean what is denoted here
  by $t \cdot A$, in conflict with the notation established earlier in
  this paper.  As a result, the last paragraph of the published
  version misstated the consequences for magnitude dimension in spaces
  with homogeneous metrics of degree smaller than $1$.}

\bibliographystyle{plain} \bibliography{pdms-corrected}

\end{document}